\RequirePackage[T1]{fontenc}
\documentclass[11pt]{article}

\usepackage{amsmath}
\usepackage{amssymb}
\usepackage[final,patch=none]{microtype}
\usepackage{graphicx}
\usepackage[hidelinks]{hyperref}
\usepackage{geometry}
\usepackage{setspace}
\usepackage{authblk}
\usepackage[switch]{lineno}
\usepackage{caption}
\usepackage{subcaption}
\usepackage{newtx}
\usepackage[utf8]{inputenc}
\DeclareUnicodeCharacter{00E9}{\'e}
\DeclareUnicodeCharacter{0308}{\"{}}

\usepackage[
  backend=biber,
  style=authoryear,
  sorting=nyt,
  giveninits=true, 
  maxbibnames=8,
  minbibnames=1,
  maxcitenames=2,
  mincitenames=1,
  sortcites=true,
  uniquename=init
]{biblatex}

\DeclareFieldFormat[article]{title}{#1}
\DeclareFieldFormat[book]{title}{#1}

\DeclareFieldFormat[article]{pages}{#1}
\DeclareFieldFormat[book]{pages}{#1}

\DeclareFieldFormat{doi}{%
  \ifhyperref
    {\href{#1}{\nolinkurl{#1}}}
    {\nolinkurl{#1}}%
}

\DeclareBibliographyDriver{article}{%
  \printnames{author}%
  \addcomma\space
  \printfield{year}%
  \adddot\space
  \printfield[title]{title}%
  \adddot\space
  \printfield{journaltitle}%
  \space
  \printfield{volume}%
  \addcomma\space
  \printfield[pages]{pages}%
  \adddot
}

\DeclareBibliographyDriver{book}{%
  \printnames{author}%
  \addcomma\space
  \printfield{year}%
  \adddot\space
  \printfield[title]{title}%
  \adddot\space
  \iffieldundef{edition}{}{\printfield{edition}\space}
  \printlist{publisher}%
  \addcomma\space
  \printlist{location}%
  \addcomma\space
  \printfield{pages}%
  \adddot
}

\DeclareBibliographyDriver{incollection}{%
  \printnames{author}%
  \addcomma\space
  \printfield{year}%
  \adddot\space
  \printfield{title}%
  \adddot\space
  In:\space
  \printfield{booktitle}%
  \adddot\space
  \iffieldundef{edition}{}{(\printfield{edition}\space ed.)}
  \printnames{editor}\space(ed[s].)
  \addcomma\space
  \printlist{publisher}%
  \addcomma\space
  \printlist{location}%
  \addcomma\space
  \printfield{pages}%
  \adddot
}

\DeclareBibliographyDriver{online}{%
  \printnames{author}%
  \addcomma\space
  \printfield{year}%
  \adddot\space
  \printfield{title}%
  \adddot\space
  \textit{Available at:}~\printfield{url}%
  \adddot\space
  \textit{Last accessed:}~\printfield{urldate}%
  \adddot
}

\AtEveryBibitem{%
  \clearfield{issn}%
  \clearfield{isbn}%
  \clearfield{url}%
  \clearfield{number}%
  \clearfield{month}%
}

\DeclareNameAlias{author}{family-given}

\usepackage{lipsum}

\begin{document}


\begin{titlepage}
\begin{flushleft}

\textbf{\LARGE{Deciphering culprits for cyanobacterial blooms and lake vulnerability in north-temperate lakes}}

\vspace{0.5cm}

\textbf{Type of article:} Research paper

\vspace{0.5cm}

\textbf{Authors:}

Jacob Serpico$^{1,a}$, B. A. Zambrano-Luna$^{1,a}$, Christopher M. Heggerud$^{2}$, Alan Hastings$^{2}$, Russell Milne$^{1,*}$, Hao Wang$^{1,*}$.

\vspace{0.5cm}

\textbf{Affiliations:}

$^1$Interdisciplinary Lab for Mathematical Ecology and Epidemiology \& Department of Mathematical and Statistical Sciences, University of Alberta, Edmonton, Alberta T6G 2R3, Canada

$^2$Department of Environmental Science and Policy, University of California - Davis, California 95616, United States

$^a$These authors contributed equally to this work.

$^*$Corresponding authors. Emails: \href{mailto:rmilne@ualberta.ca}{rmilne@ualberta.ca}, \href{mailto:hao8@ualberta.ca}{hao8@ualberta.ca}. Both of these authors supervised this work.

\vspace{0.5cm}
\textbf{Data accessibility statement:} \newline
The data that supports the findings of this study is available on GitHub with the identifier \href{https://github.com/B-A-Zambrano-Luna/Deciphering_culprits_CB/tree/main}{Deciphering culprits CB}.

\vspace{0.5cm}

\textbf{Author Contributions:} \newline
\textbf{Jacob Serpico:} Methodology, Software, Validation, Formal analysis, Investigation, Data Curation, Writing - Original Draft, Writing - Review \& Editing, Visualization, and Project administration. \textbf{B. A. Zambrano-Luna:} Methodology, Software, Validation, Formal analysis, Investigation, Data Curation, Writing - Original Draft, and Visualization. \textbf{Christopher M. Heggerud:} Validation, Writing - Original Draft, Writing - Review \& Editing. \textbf{Alan Hastings:} Validation. \textbf{Russell Milne:} Conceptualization, Methodology, Data Curation, Supervision, Writing - Original Draft, Writing - Review \& Editing, Funding acquisition. \textbf{Hao Wang:} Conceptualization, Methodology, Validation, Supervision, Funding acquisition.

\vspace{0.5cm}

\textbf{Competing Interests:}
There are no competing interests.

\vspace{0.5cm}


\textbf{Main text word count:} 5000 words.

\vspace{0.5cm}

\textbf{Figures:} Color should be used for all figures in main manuscript in print.





\end{flushleft}
\singlespacing
\begin{abstract}
Harmful cyanobacterial blooms (CBs) are increasingly prevalent worldwide, posing significant environmental and health concerns. We derive a stoichiometric model describing the population dynamics and toxicity of cyanobacteria in north-temperate freshwater ecosystems. Our model quantifies the hypoxic effects of CBs on fish mortality and evaluates the impact of microcystin-LR (MC-LR) on aquatic macro-invertebrates, phytoplankton, and fish species. Analyzing data from diverse north-temperate lakes with varying physical characteristics, we identify eutrophication as a pivotal catalyst in bloom proliferation. Under predicted warming scenarios coupled with increased eutrophication, peak MC-LR concentrations will surge dramatically, and blooms will occur earlier in the year. We uncover severe bioaccumulation of MC-LR in higher trophic species; the response to CBs among fish at intermediate trophic levels was heterogeneous across lakes. We compare our model against observations from several north-temperate lakes, demonstrating its robustness and applicability. Our insights are critical for informing targeted interventions to mitigate CBs.
\end{abstract}

\vspace{0.5cm}

\begin{flushleft}
\end{flushleft}

\end{titlepage}

\singlespacing
\twocolumn  
\normalsize
\newpage
\section*{1. Introduction}

In recent years, harmful cyanobacterial blooms (CBs) have become an increasing concern in lakes worldwide \parencite{huisman2018, taranu2012}. These blooms significantly impact aquatic ecosystems and human health due to the toxic chemicals, namely cyanotoxins, that cyanobacteria regularly produce. Cyanotoxins are known to cause adverse health effects in fish that ingest them \parencite{chorus2021}, livestock that drink contaminated water \parencite{stewart2008}, and humans that consume affected species \parencite{poste2011}. The decomposition of dead cyanobacteria following a bloom can also deplete oxygen in lakes, further exacerbating the environmental impact of CBs \parencite{paerl2013, paerl2014}. Due to the multitude of harmful effects that CBs have on lake ecosystems, recent work has been conducted to identify their causative factors and forecast them ahead of time \parencite{zhao2019, loewen2020, heggerud2024}. However, predictive models that explicitly incorporate the effects of cyanobacterial toxins on higher trophic levels and other components of lake food webs still need to be explored \parencite{soares2021, bhagowati2019, frassl2019}.

Factors that influence CBs are many and diverse. Water temperature, for instance, is highly important for determining the growth of cyanobacteria \parencite{descy2016, briand2004} because peak cyanobacterial growth rates occur at water temperatures higher than what is typically observed in nature \parencite{konopka1978}. Climate change is forecast to raise lake water temperatures, thus creating a more favourable environment for higher cyanobacterial growth rates \parencite{duan2022,desgue2023}. The growth of cyanobacteria and the prevalence of CBs are strongly correlated to nutrient availability and sunlight exposure \parencite{huisman2018, taranu2012}. Furthermore, CBs are influenced by many additional biotic and abiotic factors, such as anthropogenic nutrient pollution, water clarity, and temperature \parencite{loewen2020, giani2020}. As with temperature, increased phosphorus input from agricultural or industrial sources can exacerbate CBs \parencite{watson1997}. These factors necessitate a comprehensive understanding of how future temperature and nutrient loading patterns may impact cyanobacterial prevalence in aquatic ecosystems, especially in lakes already threatened by cyanobacteria.

CBs have a multifaceted negative influence on lake ecosystems, including producing cyanotoxins and depleting oxygen levels. Microcystins, the most common cyanotoxins, have garnered widespread attention due to their various toxic effects on a broad range of species \parencite{bouaicha2019, brophy2019, shamsollahi2018, guan2020}. \textit{Microcystis aeruginosa} is a freshwater cyanobacteria species responsible for the dominant production of microcystins \parencite{huisman2018}. Furthermore, significant microcystin-LR (MC-LR) concentrations are readily observed in lakes following CBs, highlighting the urgency of this environmental challenge \parencite{paerl2013, kotak1996, loewen2020}. The detrimental effects of MC-LR on fish health, including liver and muscle damage, have been well-documented, with implications for both aquatic ecosystems and human health \parencite{banerjee2021, malbrouck2006, shahmohamadloo2022, pearson2010}. Existing literature suggests a strong bioaccumulative effect of MC-LR and other cyanotoxins \parencite{hardy2015, kotak1996}, highlighting the threat of CBs to humans and other terrestrial or aquatic species.

\begin{figure}[ht]
\centering
\includegraphics[width=\linewidth]{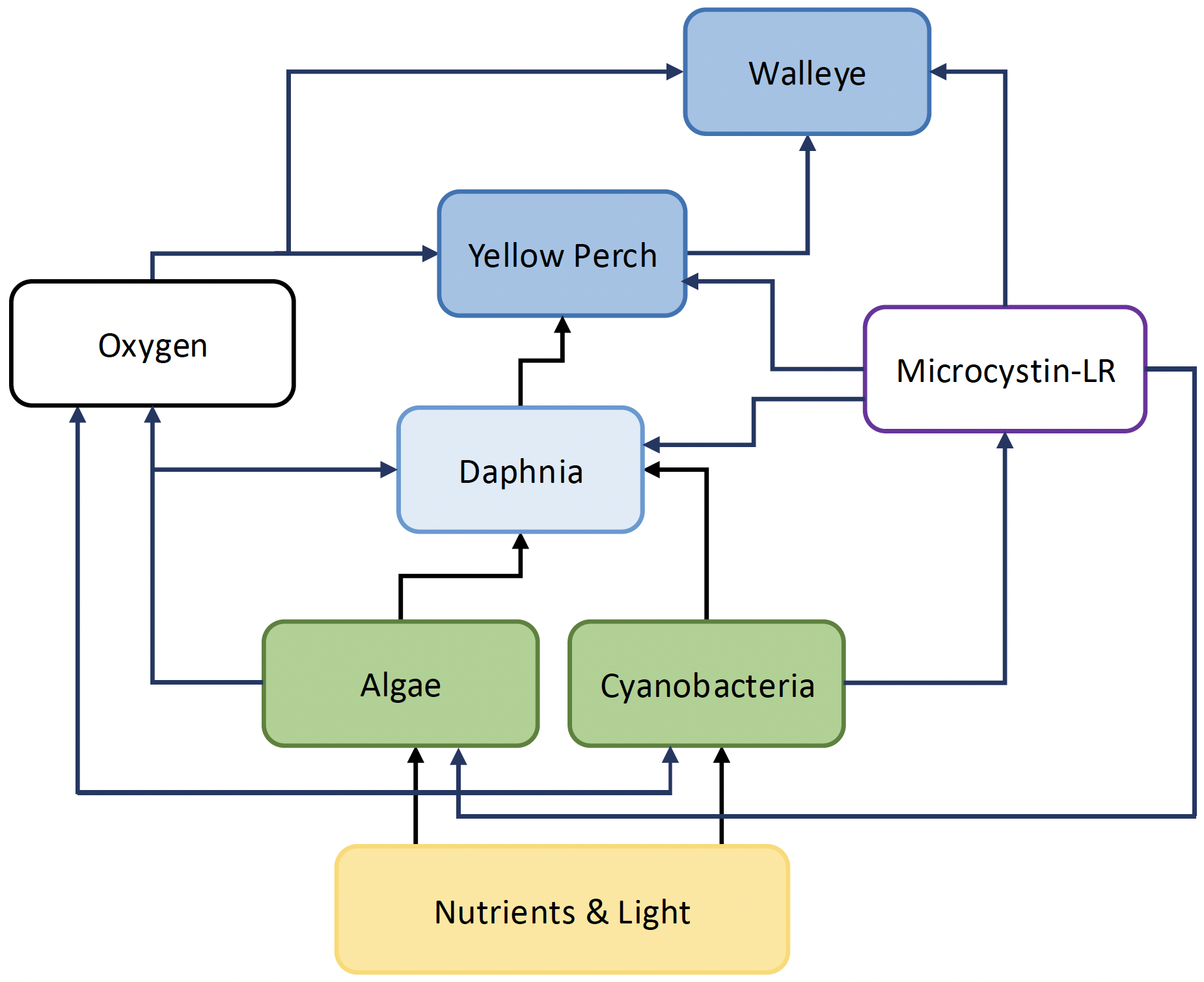}
\caption{Schematic of production and consumption within the model. Based on the principles of ecological stoichiometry, we consider the input of nutrients and light on the growth of cyanobacteria and algae. We additionally introduce predator-prey dynamics of common freshwater species in north-temperate lakes: daphnia (which consume phytoplankton), yellow perch (which consume daphnia), and walleye (which eat the smaller yellow perch). This model structure also includes the bioaccumulation of toxins produced by cyanobacteria and absorbed into the various species. Finally, we assume oxygen is produced by cyanobacteria and algae through photosynthesis and consumed by daphnia and fish species and cyanobacterial biochemical oxygen demand.}
\label{fig:foodweb}
\end{figure}

In this study, we develop a comprehensive and versatile model that considers the interactions between cyanobacteria, sunlight, and nutrients (specifically phosphorus) in a lake ecosystem; a schematic of our model is shown in Figure \ref{fig:foodweb}. We additionally consider the \textit{sander vitreus} (walleye) and \textit{perca flavescens} (yellow perch) populations, representing piscivorous (predator) and omnivorous (prey) fish species, respectively. Lastly, we consider daphnia (representing zooplankton) and algae (representing unicellular phytoplankton that do not produce toxins), the dynamics of dissolved oxygen, the production of toxins (namely MC-LR) by cyanobacteria, and the uptake of toxins by each aforementioned species. We factor in the impact of temperature on the death, metabolism, and reproduction rate of all lake species in the model, interpolating satellite data to describe water temperature and epilimnion depth. Our model is based on the principles of ecological stoichiometry, a robust modelling framework explicitly based on resource flows \parencite{sterner2002} that has been used to study algae-bacteria interactions and the impacts of toxin stress on trophic dynamics \parencite{hardy2015}. From empirical data, we simulate how temperature and phosphorus variations will affect CBs and the health of fish populations in lakes with a wide variety of physical characteristics consistent with biological ranges. We showcase the practical applicability of our model by offering actionable management strategies for various north-temperate lakes. We fit our model to data from different lakes within Canada and the United States, allowing us to estimate key components of north-temperate lake ecosystems effectively.

\section*{2. Materials and Methods}
Our study aims to predict how drastically CBs would be exacerbated by anthropogenic stressors (temperature variation due to climate change and nutrient loading), to ascertain the latent effects of these future blooms on other lake organisms, and to determine how physical lake characteristics can affect their resilience to heightened cyanobacterial activity. To answer these questions, we developed a dynamical system model based on the principles of ecological stoichiometry. The model represents a lake food web with four trophic levels, incorporates the effects of cyanobacteria on other lake organisms (secretion of MC-LR and oxygen depletion) and explicitly features ambient water temperature, phosphorus concentration, and other physical parameters relevant to cyanobacterial population dynamics (e.g. epilimnion depth, water exchange rate). The Supplementary Information, SI 1.11, contains comprehensive details on how the model was built.

\begin{figure*}[ht]
    \centering
    \includegraphics[width=0.9\linewidth]{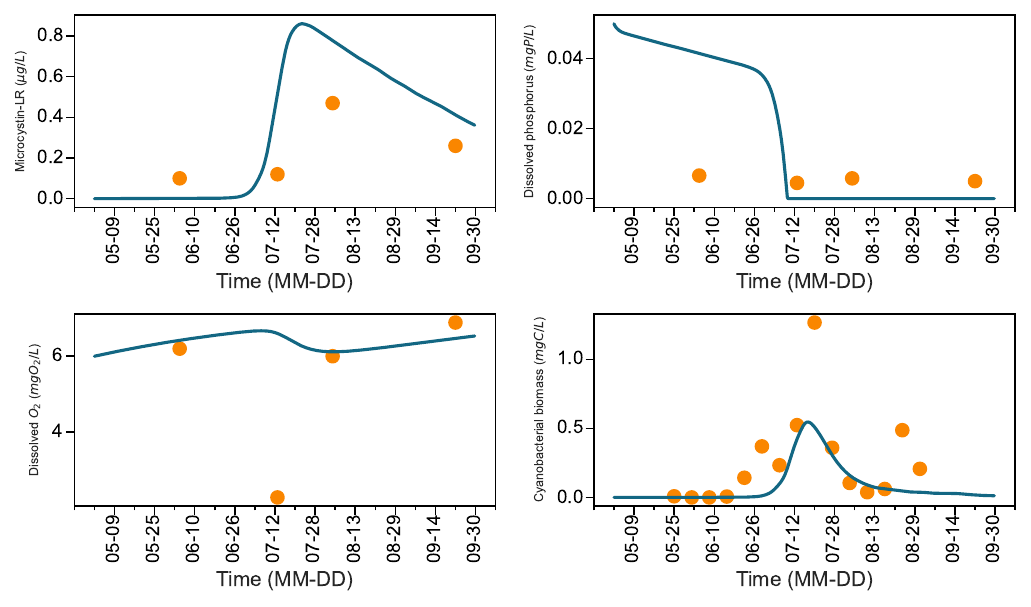}
            \label{fig:figure2}
    \caption{Predicted MC-LR, dissolved phosphorus, dissolved oxygen, and cyanobacteria (respectively in the top left, top right, bottom left, and bottom right panels) after fitting our dynamical system model to data from Pigeon lake, Alberta, Canada, 2021. The orange points represent in-situ measurements and the blue curves represent the model's solution.}
    \label{fig:my_four_plots}
\end{figure*}

\subsection*{2.1. Model formulation}

We describe the dynamics of a freshwater ecosystem, interconnecting several differential equations. For the lower trophic levels in our model (cyanobacteria and algae), we use an approach grounded in ecological stoichiometry, focusing on macronutrients, energy, and temperature interactions. We assume that phytoplankton species are grouped into either toxin producing cyanobacteria or non-toxin producing algae. Because not all unicellular phytoplankton produce toxins using these two phytoplankton gives realistic summary diet for macro-invertebrates and herbivorous fish species with respect to our overall goal of studying the toxic effects in higher trophic levels.

Further complexification of the model through more functional phytoplanktonic groups is deemed unnecessary for our purposes and unrealistic for data availability. We use a Droop form for nutrient uptake and a Monod form for light-dependent growth \parencite{morel1987, kirk2010}, both empirically supported~\parencite{wang2022}, which capture the nuanced interactions between abundances of cyanobacteria and algae, light attenuation (modelled using the Lambert-Beer law \parencite{wang2007, huisman1994}), and nutrient availability. This approach is refined to include specific adaptations of cyanobacteria, such as their optimized growth in low-oxygen conditions, by adjusting the standard growth equations to reflect these organisms' distinct ecological roles and environmental responses \parencite{heggerud2020}. Droop kinetics for nutrient uptake and a Monod form for light-dependent growth captures the inherent differences in energy conversion. Light is quickly absorbed and turned into usable energy, and the Monod form is ubiquitous for the correlation between light and growth kinetics \parencite{heggerud2020, grima1994, aiba1982, eilers1993}. Empirical evidence suggests that the Droop form describes data more accurately for nutrient uptake despite the increased mathematical complexity \parencite{wang2022}. Temperature dependence is explicitly modelled using a cardinal temperature model to describe the optimal growth temperatures for algae and cyanobacteria. This assumption is critical given the varying climatic conditions \parencite{descy2016, briand2004, vanderlinden2010} and their downstream effects on phytoplankton growth rates \parencite{cui2024}. We chose MC-LR as the representative cyanotoxin for our model due to its dominance and potence in aquatic ecosystems \parencite{tidgewell2010, NCBI2024}. In addition to cyanobacteria and algae, our model features several common species that make up a lake food web: daphnia, yellow perch, and walleye, which are represented using consumer-resource dynamics. The model also accounts for ecosystem feedback mechanisms, which are essential for understanding the resilience of aquatic species to CBs \parencite{peace2019, haney1987}. A visual representation of this model and its interactions is given in Figure \ref{fig:foodweb}.

\subsection*{2.2. Model fitting}
Numerous datasets in our data availability section were cleaned and appropriately combined for fitting. The stoichiometric model was parameterized through a rigorous process combining literature review and empirical data from three lakes in North America. As an example, we provide a model-data comparison for the fitted model to Pigeon lake data in Figure \ref{fig:my_four_plots}. Based on their dissolved phosphorus levels and data availability, we selected lakes for preliminary model fitting. Lake trophic classifications are by standards from the Canadian Water Quality Guidelines for the Protection of Aquatic Life \parencite{ccme2024}, which define oligotrophic lakes as having total phosphorus concentrations less than 10 $\mu$g/L, mesotrophic lakes between 10 and 35 $\mu$g/L, and eutrophic lakes above 35 $\mu$g/L. This diverse selection of lakes allowed us to capture a broad spectrum of ecological behaviours and interactions under varying nutrient conditions, enhancing the robustness and generalizability of our model across different lake types.

All parameters not found in the literature were estimated by minimizing the mean square error between observed and model-predicted data. See Table S.3 for all detailed parameters and their appropriate citations. We employed a differential evolution algorithm for parameter optimization, allowing for a gradient-free search within a hypercube of defined parameter constraints. Our algorithm also used a fallback strategy to ensure meaningful parameter outputs even in cases of numerical solver failure. The entire process was implemented in Python using standard packages. Further details on the fitting process and our algorithm are in the Supplementary Information; SI 1.11.

\begin{figure*}[ht]
\centering
\includegraphics[width=\textwidth]{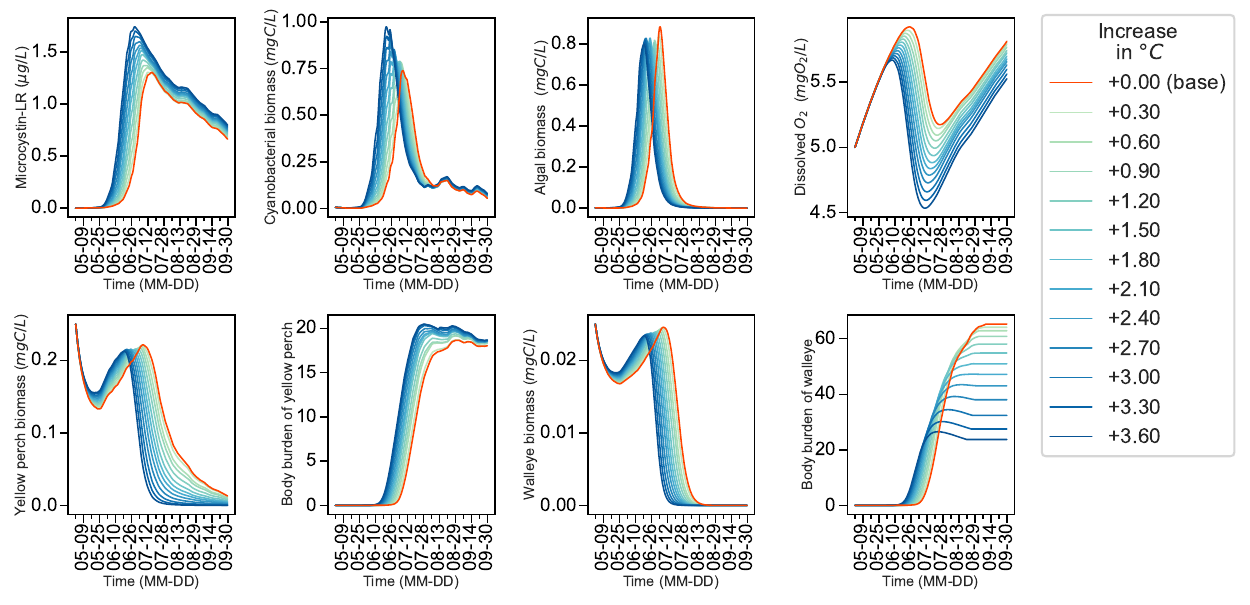}
\caption{Impacts of temperature increases concentrated during the warm season (May to September 2018) on various ecosystem metrics in Mendota lake. The plots illustrate time series of MC-LR concentration, cyanobacterial biomass, algal biomass, dissolved oxygen concentration, yellow perch biomass, and walleye biomass, and walleye and yellow perch body burden, for different temperature values compared to present conditions (+0.0$^{\circ}$C, in red).}
\label{fig:temperature}
\end{figure*}

\subsection*{2.3. Simulations}
All simulations were performed using Python, with a discretization step of $1/3$ day. We used an independent base case for each lake with parameters obtained from model fitting. A time-dependent sensitivity analysis using Sobol indices was conducted to evaluate the contributions of various physicochemical factors to cyanobacteria biomass. The factors considered included epilimnion depth, water exchange rate, turbidity, phosphorus input, and water temperature. The vulnerability index was calculated as the maximum MC-LR model output ratio under increased water temperature, epilimnion depth, and varying water exchange rates to the base case scenario (without changes). Further details are provided in the Supplementary Information.

\section*{3. Results}
We simulated dynamics for several parameter configurations of lake conditions, including variations in water temperature, phosphorus input, water exchange, and epilimnion depth. Specifically, we selected lakes with sufficient data from different mean phosphorus summertime (May 1 - Sep 30) concentrations. Referencing existing guidelines by the Canadian Council of Ministers of the Environment \parencite{ccme2024}, we chose to simulate lakes with meso-eutrophic and hypereutrophic status. These configurations reflect a diverse array of lakes in North America, and we identified critical influences on cyanobacteria concentration and MC-LR levels while considering the interactivity of algae, daphnia, yellow perch, walleye, and dissolved oxygen. 

Specifically, we simulated various lake temperatures by observing biologically reasonable ranges in north-temperate dimictic lakes. Incrementally raising the water temperature resulted in earlier and higher peaks for cyanobacteria and MC-LR in Mendota Lake in Wisconsin (Figure \ref{fig:temperature}). For Albertan lakes, Pigeon lake and Pine lake, we observed the peaks of MC-LR and cyanobacterial biomass occurring earlier. However, the height of the peaks were only increased for small to intermediate temperature increases (Figures S40, S45). These trends are mostly consistent for algal biomass, except the peaks are higher than the baseline for all increased temperatures. This subtle difference is likely because they are the competing phytoplankton species. Our simulations with increased lake temperature revealed a concerning and complex trend. They predict that height and duration of peak MC-LR concentrations vary depending on temperature. These values, reflecting a spectrum of potential scenarios influenced by population growth, economic development, and technological progress, emphasize the potential implications of future global temperature increases. 

\begin{figure*}[ht]
\centering
\includegraphics[width=1\textwidth]{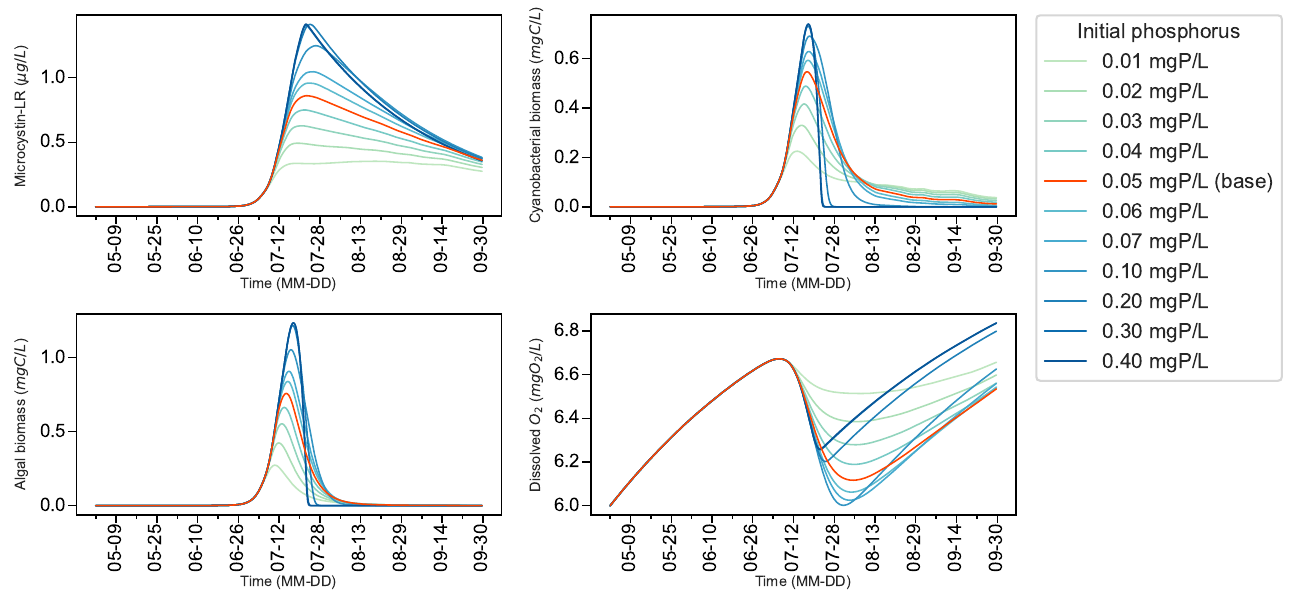}
\caption{Impacts of initial phosphorus on various ecosystem metrics in Pigeon lake. The plots illustrate time series of MC-LR concentration, cyanobacterial biomass, algal biomass, and dissolved oxygen concentration, for different initial phosphorus values compared to baseline conditions (0.05 mgP/L, in red).}
\label{fig:phosphorusnew}
\end{figure*}

Yellow perch and walleye generally follow contrasting paths as the temperature rises. Walleye biomass decreases sooner in all three lakes as temperature increases, mainly owing to the heavy body burden -- amount of toxin in an individual fish in units of mgMC-LR/mgC -- of MC-LR carried by walleye due to either higher or earlier cyanobacterial activity (Figures \ref{fig:temperature}, S40, and S45). Yellow perch biomass increases as temperature does in all lakes except Mendota (Figure \ref{fig:temperature}). As walleye are a direct predator of yellow perch, they have a significant population resurgence near the summer's end. However, the higher biomass of yellow perch is accompanied by a higher MC-LR  body burden, although this is not nearly as severe as it is in walleye due to bioaccumulation. 

Our simulation results indicate that initial phosphorus concentrations strongly influence cyanobacterial biomass and MC-LR concentrations. In all three lakes, increases in initial phosphorus dramatically increased cyanobacteria, algae and MC-LR and caused severe decreases in dissolved oxygen (Figures \ref{fig:phosphorusnew}, S35, and S43). This result delineates nutrient enrichment as a critical influence on CBs and their adverse environmental effects.

\begin{figure*}[ht]
\centering
\includegraphics[width=\textwidth]{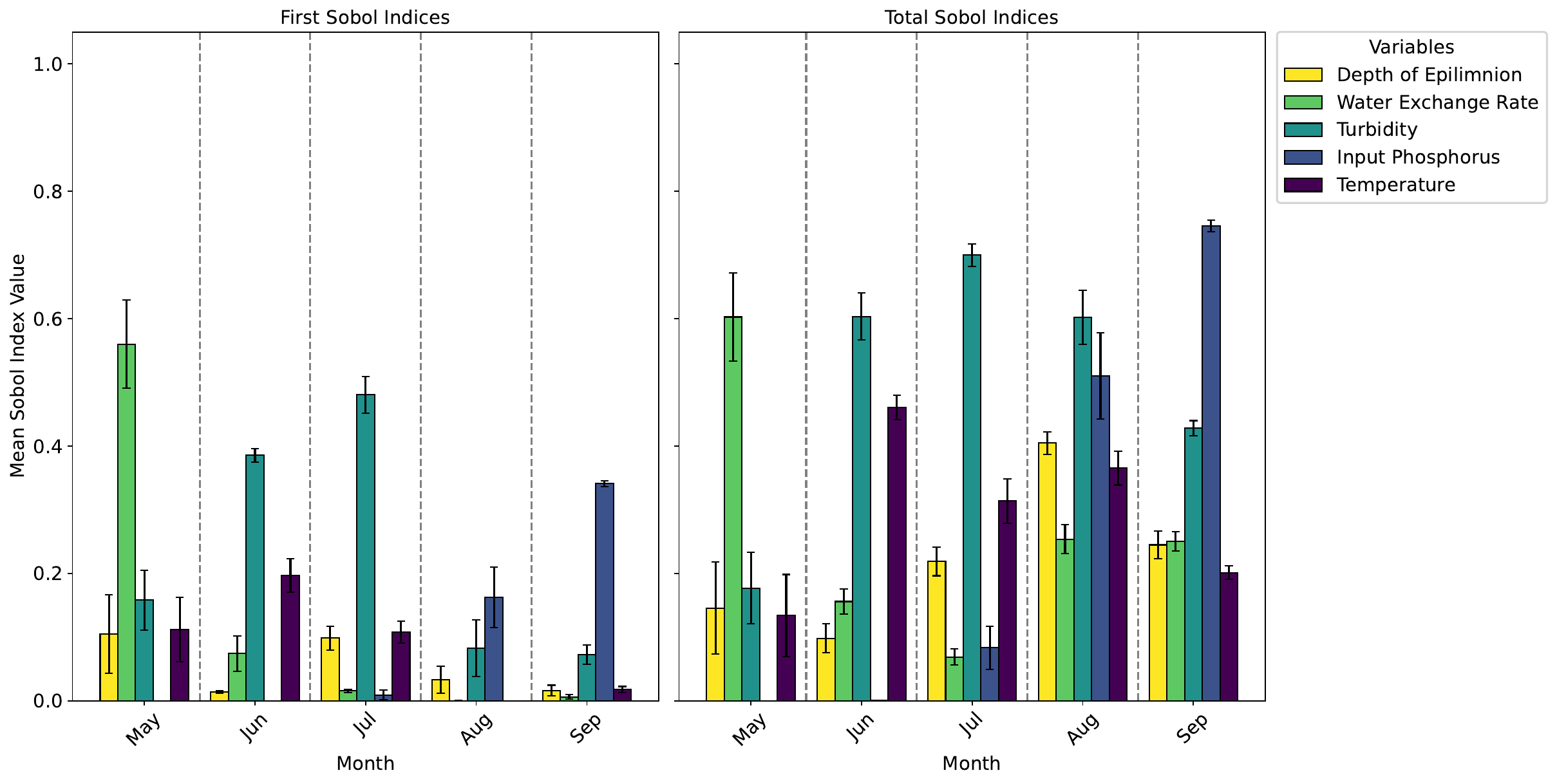}
\caption{First-order (left) and total-order (right) Sobol indices for six model features: epilimnion depth, water exchange rate, turbidity, input phosphorus, and water temperature, in Pigeon lake, 2021. These indices vary each feature and calculate their direct (first-order) and combined direct and indirect (total-order) influence on cyanobacterial biomass.}
\label{fig:CB_Sobol_index}
\end{figure*}

The influence variations on cyanobacterial biomass was also dependent on month. This is depicted in Figure \ref{fig:CB_Sobol_index}, where we calculated the first- and total-order Sobol indices. These indices serve as statistical measures, assessing the sensitivity of model outputs to changes in input parameters and quantifying how variations in specific inputs influence model outcomes \parencite{sobol2001}. In May, the water exchange rate and turbidity were the dominant contributors to cyanobacterial biomass for all tested lakes. However, from June to early August, the influence variation of epilimnion depth steadily increased, while the variation of the other features decreased. Both the first- and total-order Sobol indices of epilimnion depth peaked in the middle of August. Towards the late warm season, the contributions of water exchange rate and phosphorus began to rise, with turbidity and lake temperature following suit after the warm season ended. The individual contribution of temperature is not as substantial as other individual factors, as shown by the first-order Sobol indices. However, temperature has a strong contribution to cyanobacterial biomass when interacting with the other input parameters (Figures \ref{fig:CB_Sobol_index}, S46, and S47).

We used a vulnerability index -- the MC-LR ratio between outcomes for an increased temperature scenario and the base scenario -- to measure trends for lake characteristics most at risk of CBs. We used three increased temperature scenarios (+0.5$^{\circ}$C, +1.5$^{\circ}$C, and +3.5$^{\circ}$C) based on physically reasonable values \parencite{climatedata2023} and the results are shown in Figure \ref{fig:lakeproperties-ssp-vuln-yearlymax}. We found that lakes with low epilimnion depth are the most vulnerable to heightened MC-LR concentrations once temperatures rise, consistent with their greater propensity for CBs now \parencite{paerl2008}. The relationship between the water exchange rate and MC-LR concentration was more complex. The water exchange rate is critical for shaping bloom characteristics. However, our lower-medium warming scenarios for Pine lake had average warm-season MC-LR  concentrations rise at the same relative proportions regardless of water exchange rate if all else was equal. Further, low water exchange rates exacerbate the increases in MC-LR seen in lakes with shallow epilimnions under our severe warming scenarios, indicating that the water exchange rate will rise in importance as a factor as time goes on. 

\section*{4. Discussion}

Primary exposure pathways for cyanotoxins such as MC-LR include drinking water and recreational activities \parencite{stone2009}. The consumption of fish from contaminated water is a common but poorly studied route of human exposure \parencite{codd1999}. We show that the body burden of MC-LR was more than an order of magnitude higher in walleye than in yellow perch; these results concur with prior field studies on toxin bioaccumulation in freshwater fish \parencite{hardy2015, kotak1996}. The severe bioaccumulation observed poses critical threats to wildlife conservation and public health. This indicates that fish caught during and following CBs may pose significant health risks and be unsafe for consumption \parencite{stone2009, zhang2022}. This likewise represents a direct threat to human health through consuming contaminated fish and terrestrial animals that rely on these fish as a food source. Beyond immediate toxicity, chronic exposure to cyanotoxins could lead to long-term health effects, including liver damage \parencite{li2023}, cancer \parencite{gu2022}, and neurological disorders \parencite{metcalf2021}. Humans and other consumers of fish from lakes vulnerable to CBs are at risk of secondary poisoning, which could lead to population declines in already vulnerable species in addition to the human health issues that arise. Given these risks, comprehensive monitoring and regulation of fish harvested from waters affected by CBs is urgently needed. Equally important is the need for public health advisories and guidelines on fish consumption during and after such events, as these measures are essential to protect both human and ecological health.

\begin{figure*}[ht]
\centering
\includegraphics[width=\textwidth]{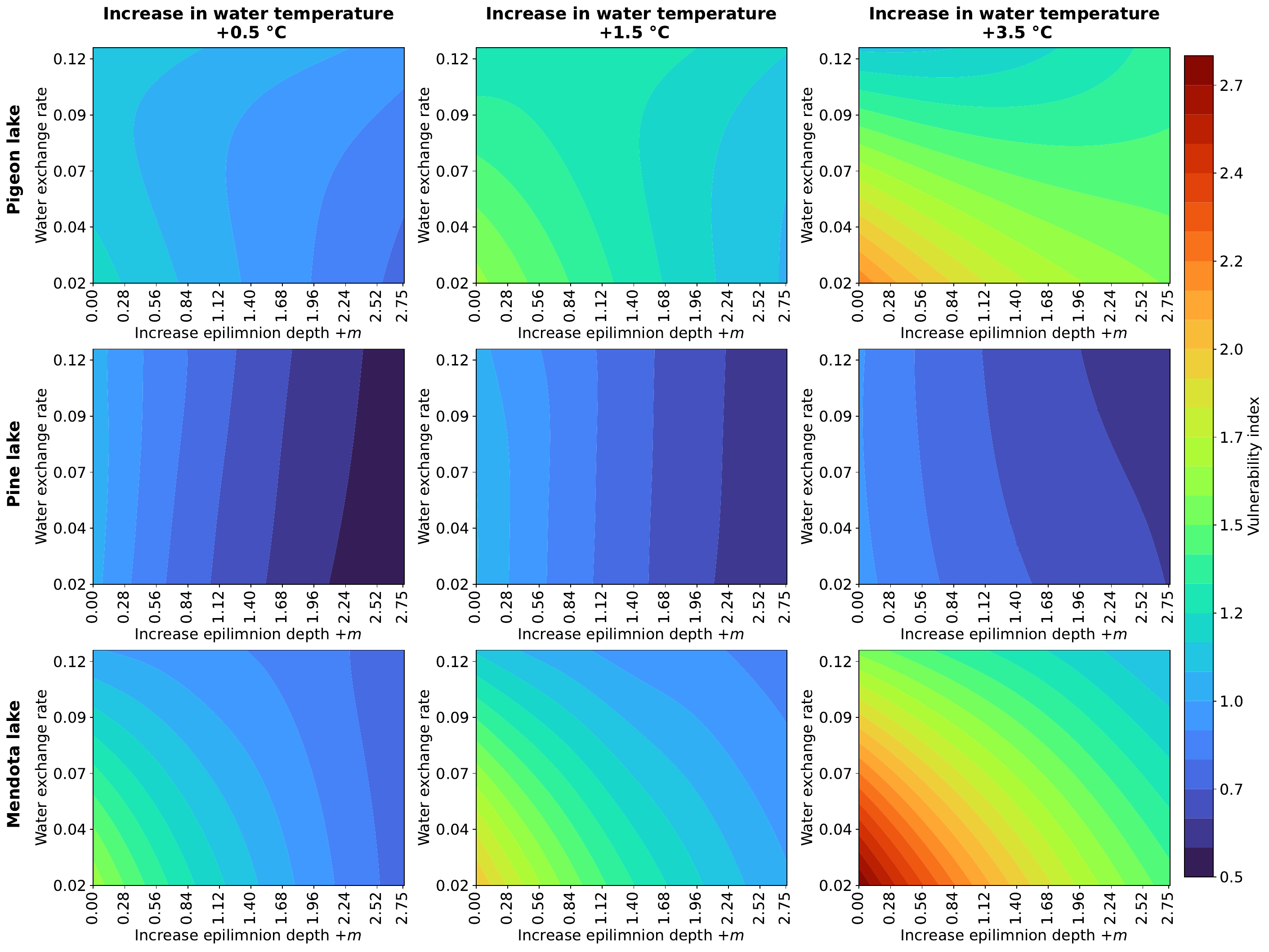}
\caption{Vulnerability indices comparing the maximum MC-LR concentration under several warming scenarios to current values. This was done for varying water exchange rates (0.02 to 0.12 m/day) and increasing the epilimnion depth (+0 to +2.7 m) for three lakes. The three rows represent Pigeon, Pine, and Mendota Lake, respectively. The columns represent increases in water temperature $+0.5$, $+1.5$, $+3.5$ $^{\circ}$C, respectively, based on physically reasonable guideline values.}
\label{fig:lakeproperties-ssp-vuln-yearlymax}
\end{figure*}

Our results delineate phosphorus and temperature as pivotal drivers in regulating the length and magnitude of CBs. The impact of temperature is particularly significant. 
Furthermore, we have shown that, in many cases, increases in temperature lead to notable increases in maximum predicted MC-LR or earlier onset of peak MC-LR concentrations in the open water season (Figures \ref{fig:temperature}, S40, and S45). Moreover, scenarios with larger temperature increases (+$1.5^{\circ}$C and +$3.5^{\circ}$C) raised MC-LR concentrations even higher (Figure \ref{fig:lakeproperties-ssp-vuln-yearlymax}). Similarly, increases in initial and input phosphorus concentrations lead to notably higher peak MC-LR concentrations (Figure \ref{fig:phosphorusnew}). This emphasizes the need for nutrient management to mitigate bloom severity and toxicity, as controlling global temperatures is comparatively futile. The growth of cyanobacteria significantly impacts the level of dissolved oxygen in lakes ~\parencite{oBoyle2016}. Our model shows that increased phosphorus concentrations can lead to significant oxygen depletion, causing strain on fish populations. This depletion is the most critical during the senescent phase of cyanobacteria, highlighting periods when lake oxygen, and thus aquatic life, are most vulnerable. We observe that the model exhibits quantitatively different behaviour from the in-situ data provided from Pine lake (Figure S25). This is likely due to the lack of spatial components in our model. In particular, our model focuses on dynamics in the well-mixed epilimnion. The measured amount of dissolved oxygen is near zero, which can not support life for extended periods of time, but could have been taken from a non-well-mixed, or deeper part of the lake. In this sense, we have less confidence that our model is capable of capturing these measurements for dissolved oxygen even though it is an important component of our model.

Several studies have explored the dynamics between environmental conditions and nutrient levels influencing CBs. We have observed similar trends that align with broader findings. For instance,  \cite{fernandez2023} expect climate change to exacerbate conditions favourable to CBs in U.S. lakes, emphasizing the complex relationship between warmer water temperatures and increased nutrient runoff. Our results agree with this relationship, showing that a temperature increase of $+3.5^{\circ}$C can cause earlier and higher peaks in MC-LR (Figures \ref{fig:temperature}, S40, S45). Moreover, \cite{li2022} found in their long-term study of Taihu lake that increased phosphorus levels were significantly correlated with the frequency and intensity of CBs. Additionally, \cite{heggerud2024} found that phosphorus concentrations were one of the critical values for predicting the onset of CBs. Our model produced similar results, where higher initial phosphorus concentrations led to substantial increases in MC-LR concentrations and peak cyanobacterial biomass. Although our results agree with a vast number of existing literature results, findings such as \cite{hellweger2022} suggest that the consideration of other macronutrients, such as nitrogen, is essential for accurately modelling cyanobacterial growth. However, we assume that phosphorus is the sole limiting nutrient for our system. This is a biologically reasonable assumption based on the water quality of north-temperate dimictic lakes \parencite{bonilla2023,wang2007}. However, the literature may benefit from models similar to ours considering nutrient co-limitation. 

Our results indicate that shallower lakes experience drastic increases in cyanobacterial levels, leading to high MC-LR levels. In contrast, lakes with deeper epilimnion layers show significantly less severe blooms, suggesting that shallow lakes may need to be a priority for managers. Figure \ref{fig:lakeproperties-ssp-vuln-yearlymax} shows shallow lakes are more vulnerable under increased temperature variations. Additionally, the water exchange rate affects bloom dynamics; a slow rate ($0.02$ m/day) maintains high toxin levels, whereas increasing the rate ($0.12$ m/day) drastically reduces MC-LR concentrations and mitigates oxygen depletion (Figure \ref{fig:phosphorusnew}). The relationship between water exchange rate and cyanobacteria has been studied previously. \cite{zhang2024} documented that an increase in the water exchange does not necessarily result in a decrease in cyanobacteria concentration, nor MC-LR, due to the quality of the exchanged water. We observed a similar effect: the relationship between the water exchange rate and peak MC-LR can become non-monotonic under warmer lake temperatures. Our simulations also indicated that lakes with low epilimnion depth have both a greater propensity for CBs now and the most significant increase in cyanobacteria activity once temperatures rise (Figure \ref{fig:lakeproperties-ssp-vuln-yearlymax}). These results suggest that the lakes currently most susceptible to CBs may become even more so, unfortunately experiencing no saturation effect relative to temperature in cyanobacterial activity.

The main factors that elicit substantial changes in bloom characteristics were identified using Sobol indices and a series of simulations for each lake. These methods were carried out with different phosphorus concentrations, temperature changes, and physical lake characteristics (Figure \ref{fig:CB_Sobol_index}). While epilimnion depth was observed to be a significant factor throughout the warm season, the water exchange rate, temperature increases, turbidity, and input phosphorus varied in both absolute and relative contributions to changes of bloom characteristics as the warm season progressed. In particular, our findings suggest that rising temperatures, specifically in May and early June (i.e. earlier summers), can promote CBs in lakes regardless of their physical characteristics. However, the effects of heightened temperatures in September (i.e. longer summers) are more muted. Instead, phosphorus input plays a more significant role during the late warm season, and more stagnant lakes with lower water exchange rates are also more prone to cyanobacteria activity during this time. These findings enable targeted interventions, such as controlling phosphorus runoff in agricultural areas or monitoring temperature changes in vulnerable lakes.

Future work should investigate the effects of other limiting factors on phytoplanktonic growth. This would include the availability of nitrogen, carbon, and micronutrients to provide a more comprehensive understanding of nutrient interactions in bloom formation. Furthermore, empirical evidence suggests that environments' spatial scale and structure can influence population interactions \parencite{cantrell2004}. It would be informative (although complex) to integrate within-lake spatial heterogeneity into our model, following in the footsteps of prior modelling work focusing on a single part of a lake ecosystem, e.g. \cite{tao2021, zhang2021}. Remote sensing using satellite imagery could provide accurate data to validate and refine such a spatial model on the scale of a lake. Additionally, as climate change continues to influence global water systems, longitudinal modelling studies that track these changes will be essential for developing better management strategies to mitigate CBs' impact and downstream effects, such as the bioaccumulation of toxins over a more extended period.

One limitation is the need for more data on oligotrophic lakes. Lakes classified as oligotrophic are typically monitored less frequently due to the current assumption that there is less propensity for CBs within them, as they intuitively have lower levels of nutrients for cyanobacteria to consume. While ample literature supports this, emerging studies show the need to actively monitor these lakes in an evolving world \parencite{reinl2021}. Lastly, as this model is already very complex, an assumption was made that walleye have no predators, only a natural death rate exacerbated by hypoxic conditions. They also do not expel the toxins back into the environment. We chose to fit only a subset of our state variables due to the availability of data. Some variables such as phosphorus and fish populations were omitted due to the lack of recent, accurate, and frequent in-situ measurements. However, all remaining parameters values, or ranges, are extracted from relevant literature, and are thus still biologically reasonable. These factors lead the model to output low walleye populations during the off-season due to their high body burden. While the shallow values depicted by the model are an artifact of its design, the broader qualitative dynamics that show a summer peak followed by a decline in the fall are as expected. These results are consistent with other milestone modelling articles \parencite{huang2015,chen2022}.


\section*{6. Code Availability}
All data, code, and supporting information used in this study is available on Zenodo (DOI: 10.5281/zenodo.15306524).

\section*{7. Acknowledgments}
Our work was primarily supported by a grant from Alberta Conservation Association. H.W. was partially supported by an NSERC Individual Discovery Grant (RGPIN-2020-03911), an NSERC Discovery Accelerator Supplement Award (RGPAS-2020-00090), and a Tier 1 Canada Research Chair. C.M.H. and A.H. were partially supported by the US NSF grant no. 2025235. C.M.H. additionally acknowledges the support of an NSERC PDF (PDF-578392-2023). We would like to express our gratitude to Peta Prokopiou for conducting a literature review during the early stage. Grammarly was used to polish writing.

\renewcommand{\bibfont}{\footnotesize}
\printbibliography

\end{document}